 
\font\smcp=cmcsc8
\font\mcp=cmcsc10


\baselineskip=14pt
\parskip=10pt
\def\Tilde{\char126\relax}

\font\eightrm=cmr8  

\font\tenrm=cmr10
\magnification=\magstephalf

\parindent=0pt
\overfullrule=0in
 
\bf
\centerline
{A 2-COLORING OF $[1,N]$ CAN HAVE $(1/22)N^2+O(N)$
MONOCHROMATIC}
\centerline
{SCHUR TRIPLES, BUT NOT LESS!}
\rm
\bigskip
\centerline{ {\mcp  Aaron Robertson and Doron Zeilberger\footnote{$^1$}
{\eightrm  \raggedright
Supported by NSF.}}}
\centerline{{Department of Mathematics, Temple University,
Philadelphia, PA 19122, USA}}
\centerline{{\tenrm aaron@math.temple.edu, zeilberg@math.temple.edu}}
\centerline{{\it Submitted:  March 3, 1998; Accepted:  March 25, 1998; Revised:  May 3, 1998}}
{}\footnote{}{\eightrm \raggedright Mathematics Classification Numbers:
Primary:  05D10, 05A16; Secondary: 04A20} 
 
{\smcp {\bf Abstract}: We prove that the minimum number (asymptotically) of
monochromatic Schur triples that a 2-coloring of $[1,n]$ can have is ${{n^2} \over {22}} + O(n)$.
This revised version fills in a minor and subtle gap discovered by M. Primak. (The
revision also corrects (at no extra cost) 
a discrepancy between the solution in the paper and the
solution obtained by Maple.  In the paper $H_{1/2}$ should be $H_0$
and $H_1$ should be $H_{1/2}$ for the solutions to agree.)}

This article is accompanied by the Maple package {\tt RON}, available from
either author's website.
 
{\bf Tianjin, June 29, 1996}: In a fascinating invited talk at the SOCA 96 
combinatorics conference organized by Bill Chen, Ron Graham proposed
(see also [GRR], p. 390):

{\bf Problem (\$100)}: Find (asymptotically) the least number
of monochromatic Schur triples $\{i,j,i+j\}$ that may
occur in a 2-coloring of the integers $1,2, \dots , n$.
 
By renaming the two colors $0$ and $1$, the above is
equivalent to the following

{\bf Discrete Calculus Problem}: Find the minimal value
of
$$
F(x_1, \dots, x_n):=
\sum_{{1 \leq i < j \leq n} \atop {i+j \leq n}}
 \left [ \,\, x_i x_j x_{i+j} + (1-x_i)(1-x_j)(1-x_{i+j}) \,\, \right ]
 ,
$$
over the $n$-dimensional
(discrete) unit cube $\{ (x_1, \dots , x_n) \vert x_i=0,1 \}$. 
We will determine {\it all} local minima (with respect to the
Hamming metric), then determine the global minimum. 
 
{\bf Partial Derivatives:}
For any function $f(x_1, \dots , x_n)$ on $\{0,1\}^n$ define
the discrete {\it partial derivatives} $\partial_r f$ by
$
\partial_r f(x_1, \dots , x_r , \dots , x_n):=
f(x_1, \dots , x_r , \dots , x_n)-
f(x_1, \dots , 1-x_r , \dots , x_n).
$
 
If $(z_1, \dots , z_n )$ is a local minimum of $F$, then 
we have the $n$ inequalities: 
$$
\partial_r F(z_1, \dots, z_n) \leq 0 \quad, \quad 1 \leq r \leq n.
$$
 
A purely routine calculation (applicable Maple routines:  {\tt diff1, dif}) shows that (below $\chi(S)$ is
1(0) if $S$ is true(false))
$$
\partial_r F(x_1, \dots , x_n)=
$$
$$
(2x_r-1) \!
\left \{ \sum_{i=1}^{n} \!x_i + \!
\sum_{i=1}^{n-r} x_i 
-(n \! - \! \left \lfloor{{r} \over {2}} \right \rfloor)-{\chi}(r\!>\! {{n} \over {2}})-
(2x_r\!-\!1)+\! x_r {\chi}(r\!> \!{{n} \over {2}}) + 1 -
(x_{{{r} \over {2}}}\!+\!x_{2r}) {\chi}(r \!\leq \!{{n} \over {2}} )\right \} .
$$
Since we are only interested in the {\it asymptotic} behavior,
we can modify $F$ by any amount that is $O(n)$. In particular,
we can replace $F(x_1, \dots , x_n)$ by 
$$
G(x_1, \dots , x_n) = F(x_1, \dots , x_n)
+ \sum_{i=1}^{n/2} x_i (x_{2i} -1) - {{1} \over {2}}\sum_{i=1}^{n} x_i.
$$
Noting that $(2x_r-1)^2 \equiv 1$ and $(2x_r-1)x_r \equiv x_r$
on $\{0,1\}^n$, we see that for $1 \leq r \leq n$,
$$
\partial_r G(x_1, \dots , x_n)=
(2x_r-1) 
\left \{ \sum_{i=1}^{n} x_i +
\sum_{i=1}^{n-r} x_i 
-(n- \left \lfloor {{r} \over {2}} \right \rfloor) - {{1} \over {2}} {\chi}(r \leq n/2) 
\right \} -{{1} \over {2}} {\chi}(r \leq n/2) - 1/2.
$$
Let $k=\sum_{i=1}^{n} x_i$.
Since at a local minimum $(z_1, \dots, z_n)$ we have 
$\partial_r G(z_1, \dots, z_n) \leq0$,
it follows that any local minimum $(z_1, \dots, z_n)$ satisfies the
 
{\bf Ping-Pong Recurrence}:
Choose $a,b \in \{0,1\}$ arbitrarily each time $\widehat{H}$ or 
$\widetilde{H}$ is used, where $\widehat{H}$ and
$\widetilde{H}$ are the following functions:
$$
\widehat{H}(y):=\cases{0,& if $y > 1/2$;\cr
             1,& if $y < 0$;\cr
             a,& if $0 \leq y \leq 1/2$.\cr}
$$
$$
\widetilde{H}(y):=\cases{0,&if $ y > 1$;\cr
		  1,&if $y < -1$;\cr
		  b,&if $-1 \leq y \leq 1$.\cr}
$$
Then we must have, for $r=n, n-1 , \dots, n-\lfloor n/2 \rfloor + 1$,
$$
z_{r}=\widehat{H}\left(k-n+ \left \lfloor {{r} \over {2}} \right \rfloor + \sum_{j=1}^{n-r} z_j\right) ,
\eqno(Right\,Volley)
$$
$$
z_{n-r+1}=\widetilde{H}\left(2k-n-1/2+ \left \lfloor {{n-r+1} \over {2}} \right \rfloor- \sum_{j=r}^{n} z_j\right) ,
\eqno(Left \,Volley)
$$
and if $n$ is odd then 
$z_{(n+1)/2}=\widehat{H}(k-n+\lfloor {{n+1} \over {4}} \rfloor + \sum_{j=1}^{(n-1)/2} z_j)$.

These equations determine a solution (depending upon the choices of
the $a$'s and $b$'s made along the way), $z$ (if it exists), 
in the order $z_n,z_1,z_{n-1},z_2, \dots $. When we solve
the Ping-Pong recurrence we forget the fact that
$\sum_{i=1}^{n} z_i=k$. Most of the time a solution
will not satisfy this last condition, but when it does, we
have a genuine local minimum. Note that {\it any} local minimum
must show up in this way.
 
{\bf Solutions of the Ping-Pong Recurrence:}
By playing around with the Maple routine {\tt ptor2} in our
Maple package {\tt RON}
(available from either author's website), we were
able to find the following solutions, 
{\it for n sufficiently large}, to
the Ping-Pong recurrence.
As usual, for any word (or letter) $W$,
$W^m$ means `$W$ repeated $m$ times'.
 
Let $w=2k-n$, $k \neq n/2$ (this case must be dealt with separately).
By symmetry we may assume that $k\geq n/2$.
Then $0 <w \leq n$. If $w \geq n/2$ then
the only solution is $0^n$. If $w< n/2$, then
let $s$ be the unique integer $0 \leq s < \infty$, that satisfies
$n/(12s+14) \leq w < n/(12s+2)$. 
 
{\bf Case I:} If $n/8 \leq w < n/2$ then the
solutions are:
$$
0^{\lfloor{{n} \over {2}} \rfloor }1^{n-\lfloor{{n} \over {2}} \rfloor-w-c_1}0^{w+c_1}
$$
where $c_1 \in \{-1,0,1\}$.

{\bf Case II:} If $n/(12s+8) \leq w < n/(12s+2)$ then the
solutions are
$$
\cases{
0^{4w+c_1} 1^{\lfloor n/2 \rfloor - 4w - c_1} 0^{n-\lfloor n/2 \rfloor -7w-(c_2+c_3+c_4)} 
1^{6w+c_3} 0^{w+c_4}
& for $s=1$;\cr
0^{4w+c_4} (1^{6w+c_5^{s_i}} 0^{6w+c_6^{s_i}})^{s/2 } Q
(0^{6w+c_7^{s_i}} 1^{6w+c_8^{s_i}})^{s/2} 0^{w+c_9}
& for $s > 1$.\cr}
$$
where the $c_j$'s and $c_j^{s_i}$'s are bounded constants (independent of $n$)
and $Q$ can be
an (almost) arbitrary mix of $r$ zeroes and ones 
(where $r$ is the unique integer such that the length of this interval is $n$).
Further, the number of ones in $Q$ is 
at most $12w$.  Notation: (1) the $c_j^{s_i}$'s
can be different constants with $i$ ranging from $1$ to $s/2$;
(2) if $s$ is odd $(a b)^{s/2}$ is
$(a b)^{(s-1)/2} a$.

{\bf Case III:} If $n/(12s+14) \leq w < n/(12s+8)$ then the
solutions are
$$
\cases{
0^{4w+d_1} 1^{n - 5w - (d_1+d_2)} 0^{w+d_2}
& for $s=0$;\cr
0^{4w+d_3} (1^{6w+d_4^{s_i}} 0^{6w+d_5^{s_i}})^{ s/2 } Q
(0^{6w+d_6^{s_i}} 1^{6w+d_7^{s_i}})^{s/2 } 0^{w+d_8}
& for $s > 0$.\cr}
$$
where the $d_j$'s and $d_j^{s_i}$'sare bounded constants (independent of $n$)
and $Q$ can be an (almost) arbitrary mix of $r$ zeroes and ones, with the
number of ones in $Q$ at most $6w$. 

{\bf Case IV:} if $w=0$ (i.e. $s=\infty$),
the solutions are:
$$
0^{g_1} (1^{g_2^{n_i}} 0^{g_3^{n_i}})^{n/(2G_1)} Q 
(0^{g_4^{n_i}} 1^{g_5^{n_i}})^{n/(2G_2)}
$$ 
where $g_1 \in \{0,1,2\}$, the other $g_i$'s and $g_i^{n_i}$'s
are bounded between $3$ and $11$,
$Q$ is an (almost) arbitrary mix of $r$ zeroes and ones
with the number of ones bounded between $0$ and $22$, $G_1 = \sum_{i} (g_2^{n_i}
+g_3^{n_i})$, and $G_2 = \sum_{i} (g_4^{n_i}+g_5^{n_i})$.

{\bf Proof:} Routine verification!

Now it is time to impose the extra condition that
$\sum_{i=1}^{n} z_i=k$ ($=(w+n)/2)$. With Cases I and II
a routine calculation yields a contradiction of the applicable range of $w$
when $n$ is sufficiently large.
For Case III, a routine calculation yields a local minimum of
$w=n/11$ if $s=0$.  If $s>0$ argue as follows.  Let $t$ be the number of
$1$'s in $Q$.  Recall that $r$ is the total number of $0$'s and $1$'s in
$Q$.  Let $w_c(s) = n/(12s+c)$ where we must have $8 \leq c \leq 14$.  
Since we need $\sum_{i=1}^{n} z_i=k$ ($=(w+n)/2)$, we see that
$6w_c(s)s+t= n(12s+c+1)/(24s+2c)$ gives $t=(c+1)w_c(s)/2$.  
Further, since
the number of $1$'s in $Q$ is bounded by $6w_c(s)$, we find that we must have
$c \leq 11$.  We also must have $r=n - w_c(s)(12s+5)$, by the definition of $r$.
Using the simple
inequality $r \geq t$, we have $n - w_c(s)(12s+5) \geq (c+1)w_c(s)/2$.  From
this deduce that $c \geq 11$.  Hence we must have $c=11$ at a local minimum.
Thus the local minimums for Case III, $s>0$, are $w_s=n/(12s+11)$.
Case IV gives infinitely many local minimums.
Hence
 
{\bf The Local Minima Are Asymptotically Equivalent (mod $O(n)$) to:}
$$
\cases{
Z_s:=0^{4w_s} (1^{6w_s} 0^{6w_s})^{{{s} \over {2}}} 1^{6w_s} 
(0^{6w_s} 1^{6w_s})^{{{s} \over {2}}}0^{w_s}
& for $0 \leq s < \infty$ (where $w_s:={{n} \over {12s+11}}$),\cr
Z_{\infty}^t= (0^t 1^t)^{n/(2t)}
& for $3 \leq t \leq 11$ \cr}
$$
A routine calculation [R] shows that for $0 \leq s < \infty$
$$ 
F(Z_s)={{12s+8} \over {16(12s+11)}}n^2 +O(n),
$$
which is strictly increasing in $s$.  An easy 
calculation shows $F(Z_\infty^t)=
(1/16)n^{2} + O(n)$ for any natural number $t$.
 
{\bf ...And The Winner Is:}\quad
$
Z_0=0^{4n/11}1^{6n/11}0^{n/11} \,\,
$
setting the world-record of $(1/22)n^2+O(n)$.
 
{\bf An Extension}
Here we note that our result implies a good upper bound for the 
general r-coloring of the first $n$ integers;  if we r-color the
integers (with colors $C_1 \dots C_r$) from 1 to $n$ then
the minimum number of monochromatic Schur triples is bounded above by
$$
{{n^2} \over {2^{2r-3}11}} + O(n).
$$
This comes from the following coloring:
$$
\cases{
Color(i)=C_j
&if ${{n} \over {2^j}} < i \leq {{n} \over {2^{j-1}}}$ \quad for $1 \leq j \leq r-2$,\cr
Color(i)=C_{r-1}
&if $1 \leq i \leq {{4n} \over {2^{r-2}11}}$ or
   $ {{10n} \over {2^{r-2}11}} < i \leq {{n} \over {2^{r-2}}}$,\cr
Color(i)=C_r
&if ${{4n} \over {2^{r-2}11}} < i \leq {{10n} \over {2^{r-2}11}}$.\cr}
$$

{\bf Note:} Tomasz Schoen[S],
a student of Tomasz Luczak, has independently solved this problem.
 
{\bf REFERENCES}
 
[GRR] R. Graham, V. R$\ddot{o}$dl, and A. Rucinski,
{\it On Schur properties of random subsets of integers},
J. Number Theory {\bf 61} (1996), 388-408.

[RZ] A. Robertson,
{\it On the asymptotic behavior of Schur triples},
[{\tt www.math.temple.edu/\Tilde aaron}].
 
[S] T. Schoen, {\it On the number of monochromatic Schur triples}, in 
preparation, 
[{\tt wtguest3@informatik.uni-kiel.de}].

\bye